\documentclass{elsart}

\journal{Journal of Number Theory}

\usepackage{amsmath}
\usepackage{graphicx}

\begin{document}

\newcommand {\Res}{\mathop{\textup{Res}}}
\newcommand {\Li}{\mathop{\textup{Li}}\nolimits}
\newtheorem {Theorem}{Theorem}

\begin{frontmatter}

 \title{On the integral of the fourth Jacobi theta function}
 \author{Istv\'an Mez\H{o}\thanksref{a}}
 \thanks[a]{Present address: University of Debrecen, H-4010, Debrecen, P.O. Box 12, Hungary}
 \address{Department of Applied Mathematics and Probability Theory, Faculty of Informatics, University of Debrecen, Hungary}
 \ead{mezo.istvan@inf.unideb.hu}
 \ead[url]{http://www.inf.unideb.hu/valseg/dolgozok/mezoistvan/mezoistvan.html}

\begin{abstract}We generalize the Raabe-formula to the $q$-loggamma function. As a consequence, we get that the integral of the logarithm of the fourth Jacobi theta function between its least imaginary zeros is connected to the partition function and the Riemann zeta function.
\end{abstract}

\end{frontmatter}

\section{Introduction}

The main result of the paper is the next integral formula for the fourth Jacobi theta function:
\[\int_{-x^*}^{x^*}\log\vartheta_4(x,q)dx=i\left[\zeta(2)-\log q\log\sum_{n=0}^\infty P(n)q^{2n}\right].\]
Here $x^*=\frac12i\log q$ is the least zero (on the imaginary axis) for the Jacobi theta function
\[\vartheta_4(x,q)=\sum_{n=-\infty}^\infty(-1)^nq^{n^2}e^{2nix},\]
$P(n)$ is the partition function of the natural numbers \cite{Andrews}, $0<q<1$ and $i=\sqrt{-1}$.

\section{Preliminaries}

\subsection{Raabe's formula}

In 1840 J. L. Raabe \cite{Raabe} proved that for the Euler $\Gamma$ function
\[\int_0^1\log\Gamma(x+t)dx=\log\sqrt{2\pi}+(t\log t-t)\quad(t\ge 0).\]
This implies the special case
\[\int_0^1\log\Gamma(x)dx=\log\sqrt{2\pi},\]
and an immediate consequence is that
\[\int_0^1\log\Gamma(x)\Gamma(1-x)dx=\log2\pi.\]
(See \cite{Amdeberhan} for an elementary proof of this special case.)
We shall prove the appropriate integral formula for the Jackson $q$-gamma function. Then we show that the Jacobi triple product identity connects the $q$-gamma function to $\vartheta_4$ and our main formula will follow.

\subsection{The Jackson's $q$-gamma function}

F. H. Jackson defined the $q$-analogue of the standard Euler $\Gamma(z)$ function as \cite{Jackson,Moak}
\begin{equation}
\Gamma_q(z)=\frac{(q^{-1};q^{-1})_\infty}{(q^{-z};q)_\infty}(q-1)^{1-z}\quad(q>1)\label{gdammadef}
\end{equation}
with $(x;q)_\infty=(1-x)(1-qx)(1-q^2x)\cdots$.

\subsection{The zeta regularized product}

Let us consider a sequence $\textbf{a}=(a_1,a_2,\dots)$. Its zeta regularized product is denoted and defined by \cite{Kurokawa,Kurokawa2}
\[\widehat{\prod_{n=1}^\infty}a_n=\exp(-\zeta'_{\textbf{a}}(0)).\]
Here
\[\zeta'_{\textbf{a}}(s)=\sum_{n=1}^\infty a_n^{-s}\]
is the zeta function associated with the sequence $\textbf{a}$. It is assumed that $\zeta_{\textbf{a}}(s)$ has an analytic
continuation to a region containing $s=0$ and further that it is holomorphic at this point.

M. Lerch \cite{Lerch} proved that for $\textbf{a}=(1+x,2+x,\dots)$
\[\widehat{\prod_{n=1}^\infty}(n+x)=\frac{\sqrt{2\pi}}{\Gamma(x)}\quad(x>0).\]
This comes from the fact that the sequence $\textbf{a}$ above has the associated zeta function
\[\zeta_{\textbf{a}}(s)=\zeta(s,x)=\sum_{n=1}^\infty (n+x)^{-s}\]
(which is the well known Hurwitz zeta function) and that
\[\zeta'(0,x)=\log\frac{\Gamma(x)}{\sqrt{2\pi}}.\]

Our considerations need a more general form of the zeta regularization. This was worked out in 2005 by N. Kurokawa and N. Wakayama \cite{Kurokawa}. They assumed that the $\zeta_{\textbf{a}}(s)$ function is meromorphic at $s=0$ with the Laurent expansion
\[\zeta_{\textbf{a}}(s)=\sum_{m}c_m(\textbf{a})\textbf{s}^m.\]
Then the generalized zeta regularized product is defined by
\begin{equation}
\widehat{\prod_{n=1}^\infty}a_n=\exp(-c_1(\textbf{a}))=\exp\left(-\Res_{s=0}\frac{\zeta_{\textbf{a}}(s)}{s^2}\right).\label{GZRP}
\end{equation}
See \cite{Mizuno} for nice applications.

We introduce the short and standard notation
\[[n]_q=\frac{q^n-1}{q-1}.\]
With this abbreviation the zeta function associated with the sequence $\textbf{a}=([x],[1+x]_q,[2+x]_q,\dots)$ is
\[\zeta_{\textbf{a}}(s)=\zeta_q(s,x)=\sum_{n=0}^\infty [n+x]_q^{-s}.\]

Then the generalized zeta regularized product reads as
\begin{equation}
\widehat{\prod_{n=0}^\infty}[n+x]_q=\widehat{\prod_{n=1}^\infty}\frac{q^{n+x}-1}{q-1}=\frac{C_q}{\Gamma_q(x)}\quad(q>1),\label{NK}
\end{equation}
where
\begin{equation}
C_q=q^{-\frac{1}{12}}(q-1)^{\frac12-\frac{\log(q-1)}{2\log q}}(q;q)_\infty.\label{Cqdef}
\end{equation}
This is the second theorem of Kurokawa and Nakayama in \cite{Kurokawa} and this will be our main tool. (A more general form of this theorem is presented in \cite{Mizuno}.)

We split the presentation to two sections. The next one contains the generalized Raabe's formula, the other contains the proof of the integral formula of $\vartheta_4$.

\section{Integral of the $q$-loggamma function}

The $q$-analogue of Raabe's theorem is given:

\begin{Theorem}For any $t>0$ and $q>1$
\begin{equation}
\int_0^1\log\Gamma_q(x+t)dx=\label{Thm1}
\end{equation}
\[\log C_q-\frac{1}{2q^t\log q}\left[\frac{1-q^t}{1-q^{-t}}(2\Li_2(q^{-t})+\log^2(1-q^{-t}))+\right.\]
\[\left.2\frac{1-q^t}{1-q^{-t}}\log\frac{1-q}{1-q^t}\log(1-q^{-t})-q^t\log^2\frac{1-q}{1-q^t}\right].\]
In special,
\begin{equation}
\int_0^1\log\Gamma_q(x)dx=\frac{\zeta(2)}{\log q}+\log\sqrt{\frac{q-1}{\sqrt[6]{q}}}+\log(q^{-1};q^{-1})_\infty.\label{Thm1spec}
\end{equation}
\end{Theorem}
Here $\Li_2(z)$ is the dilogarithm function:
\[\Li_2(z)=\sum_{n=1}^\infty\frac{z^n}{n^2}.\]

To prove this theorem, we need the following Theorem, which is interesting itself.
\begin{Theorem}For all $q>1$,
\[\int_0^1\zeta_q(s,x+t)dx=\frac{(q-1)^s}{s\log q}\frac{(q^t-1)^{1-s}}{q^t}{_2}F_1(1,1;s+1;q^{-t}).\]
\end{Theorem}
Here
\[{_2}F_1(a,b;c;z)=\sum_{n=0}^\infty\frac{(a)_n(b)_n}{(c)_n}\frac{z^n}{n!}\]
is the hypergeometric function and $(a)_n=a(a+1)\cdots(a+n-1)$ is the Pochhammer symbol.

\textit{Proof of Theorem 2.}
\[\int_0^1\zeta_q(s,x+t)dx=\int_0^1\sum_{n=0}^\infty[n+x+t]_q^{-s}dx=\]
\[(q-1)^s\sum_{n=0}^\infty\int_0^1(q^{n+x+t}-1)^{-s}dx.\]
This latter integral is computed by Wolfram Mathematica:
\[\int_0^1(q^{n+x+t}-1)^{-s}dx=\frac{1}{s\log q}\left[\frac{(q^{n+t}-1)^{1-s}}{q^{n+t}}{_2}F_1(1,1;s+1;q^{-n-t})-\right.\]
\[\left.\frac{(q^{n+t+1}-1)^{1-s}}{q^{n+t+1}}{_2}F_1(1,1;s+1;q^{-n-t-1})\right].\]
Since
\[{_2}F_1(1,1;s+1;q^{-n-t})=\sum_{k=0}^\infty\frac{k!}{(s+1)_k}\frac{1}{(q^{n+t})^k},\]
\[\sum_{n=0}^\infty\int_0^1(q^{n+x+t}-1)^{-s}dx=\]
\[\frac{1}{s\log q}\sum_{n=0}^\infty\sum_{k=0}^\infty\frac{k!}{(s+1)_k}\left(\frac{(q^{n+t}-1)^{1-s}}{(q^{n+t})^{k+1}}-\frac{(q^{n+t+1}-1)^{1-s}}{(q^{n+t+1})^{k+1}}\right).\]
If we interchange the order of the summation, we see that the sum over $n$ is telescopic, so the only one term which is not cancels belongs to $n=0$. Thus the above expression simplifies to
\[\frac{1}{s\log q}\sum_{k=0}^\infty\frac{k!}{(s+1)_k}\frac{(q^t-1)^{1-s}}{(q^t)^{k+1}}=\frac{(q^t-1)^{1-s}}{sq^t\log q}\sum_{k=0}^\infty\frac{k!}{(s+1)_k}\frac{1}{(q^t)^k}.\]
This latter sum is again hypergeometric with parameters $(1,1;s+1;q^{-t})$ and we get our Theorem.

\textit{Proof of Theorem 1.} \eqref{GZRP} and \eqref{NK} together gives that
\[\int_0^1\log\Gamma_q(x+t)dx=\log C_q+\int_0^1\Res_{s=0}\frac{\zeta(s,x+t)}{s^2}dx.\]
Since the residue is taken with respect to $s$, we can carry out before the integral. Hence, by Theorem 2,
\[\int_0^1\log\Gamma_q(x+t)dx=\log C_q+\Res_{s=0}\frac{(q-1)^s}{s^3\log q}\frac{(q^t-1)^{1-s}}{q^t}{_2}F_1(1,1;s+1;q^{-t}).\]
The residue can be calculated with Mathematica. It equals to
\begin{equation}
\frac{-1}{2q^t\log q}\left[(1-q^t)\frac{\partial^2}{\partial s^2}{_2}F_1(1,1;s;q^{-t})_{s=1}+\right.\label{Res}
\end{equation}
\[\left.2(1-q^t)\log\frac{1-q}{1-q^t}\frac{\partial}{\partial s}{_2}F_1(1,1;s;q^{-t})_{s=1}-q^t\log^2\frac{1-q}{1-q^t}\right].\]
Now we deal with the partial derivatives. Simbolically,
\begin{equation}
\frac{\partial^n}{\partial s^n}{_2}F_1(1,1;s;z)=\sum_{n=0}^\infty(a)_n(b)_n\frac{\partial^n}{\partial s^n}\frac{1}{(c)_n}\frac{z^n}{n!}.\label{hypder}
\end{equation}
The Pochhammer symbol can be rewritten with the $\Gamma$ function:
\[(c)_n=\frac{\Gamma(c+n)}{\Gamma(c)},\]
whence
\begin{equation}
\frac{\partial}{\partial s}\frac{1}{(c)_n}=\frac{-1}{(c)_n}(\psi(c+n)-\psi(c)),\label{pardev1}
\end{equation}
and
\[\frac{\partial^2}{\partial s^2}\frac{1}{(c)_n}=\frac{(\psi(c+n)-\psi(c))^2}{(c)_n}-\frac{\psi'(c+n)-\psi'(c)}{(c)_n}.\]
Here
\[\psi(z)=\frac{\Gamma'(z)}{\Gamma(z)}\]
is the digamma function. When $n$ is a positive integer, then \cite[p. 13]{AndrewsAskeyRoy}
\begin{equation}
\psi(n)=\frac11+\frac12+\cdots+\frac{1}{n-1}=H_{n-1}-\gamma,\label{specval1}
\end{equation}
and 
\[\psi'(n)=-\frac{1}{1^2}-\frac{1}{2^2}-\cdots-\frac{1}{(n-1)^2}+\zeta(2)=-H_{n-1,2}+\zeta(2).\]
($H_n$ and $H_{n,2}$ are the harmonic and second order harmonic numbers, respectively. $H_0=H_{0,2}=0$.)
Now \eqref{hypder}, \eqref{pardev1} and \eqref{specval1} gives that
\[\frac{\partial}{\partial s}{_2}F_1(1,1;s;z)_{s=1}=\sum_{n=0}^\infty n!n!\frac{-H_n}{n!}\frac{z^n}{n!}=-\sum_{n=0}^\infty H_nz^n=\frac{\log(1-z)}{1-z}.\]
Similarly, for the second order derivative
\begin{equation}
\frac{\partial^2}{\partial s^2}{_2}F_1(1,1;s;z)_{s=1}=\sum_{n=0}^\infty n!n!\left(\frac{H_n^2}{n!}+\frac{H_{n,2}}{n!}\right)\frac{z^n}{n!}=\sum_{n=1}^\infty H_n^2z^n+\sum_{n=1}^\infty H_{n,2}z^n.\label{secder}
\end{equation}
By Cauchy's product, the latter sum is simply
\[\frac{1}{1-z}\sum_{n=1}^\infty\frac{z^n}{n^2}=\frac{\Li_2(z)}{1-z}.\]
The first sum can be determined easily. Note that
\[H_{n-1}^2=\left(H_n-\frac1n\right)^2=H_n^2+\frac{1}{n^2}-2\frac{H_n}{n},\]
whence
\begin{equation}
\sum_{n=1}^\infty H_{n-1}^2z^n=\sum_{n=1}^\infty H_n^2z^n+\sum_{n=1}^\infty\frac{z^n}{n^2}-2\sum_{n=1}^\infty\frac{H_n}{n}z^n.\label{tmps1}
\end{equation}
The last sum equals to \cite{BB}
\begin{equation}
\sum_{n=1}^\infty\frac{H_n}{n}z^n=\Li_2(z)+\frac12\log^2(1-z).\label{tmps2}
\end{equation}
If we temporarily introduce the function
\[f(z)=\sum_{n=1}^\infty H_n^2z^n,\]
then \eqref{tmps1} and \eqref{tmps2} implies that
\[zf(z)=f(z)+\Li_2(z)-2(\Li_2(z)+\frac12\log^2(1-z)),\]
hence
\[f(z)=\sum_{n=1}^\infty H_n^2z^n=\frac{\Li_2(z)+\log^2(1-z)}{1-z}.\]
Altogether, \eqref{secder} becomes
\[\frac{\partial^2}{\partial s^2}{_2}F_1(1,1;s;z)_{s=1}=\frac{2\Li_2(z)+\log^2(1-z)}{1-z}.\]
The partial derivatives in \eqref{Res} are determined and the first part of the Theorem is proved.

If $t\to 0$, our expression under \eqref{Thm1} simplifies:
\[\lim_{t\to 0}\frac{-1}{2q^t\log q}\left[\frac{1-q^t}{1-q^{-t}}(2\Li_2(q^{-t})+\log^2(1-q^{-t}))+\right.\]
\[\left.2\frac{1-q^t}{1-q^{-t}}\log\frac{1-q}{1-q^t}\log(1-q^{-t})-q^t\log^2\frac{1-q}{1-q^t}\right]=\]
\[\frac{-1}{2\log q}\left[-2\Li_2(1)-\lim_{t\to 0}\left(+\log^2(1-q^{-t})+2\log\frac{1-q}{1-q^t}\log(1-q^{-t})+\log^2\frac{1-q}{1-q^t}\right)\right]=\]
\[\frac{-1}{2\log q}\left[-2\Li_2(1)-\lim_{t\to 0}\left(\log(1-q^{-t})+\log\frac{1-q}{1-q^t}\right)^2\right]=\]
\[\frac{-1}{2\log q}\left[-2\Li_2(1)-\lim_{t\to 0}\log^2\left((1-q)\frac{1-q^{-t}}{1-q^t}\right)\right]=\frac{1}{2\log q}(2\zeta(2)+\log^2(q-1)).\]
Thus \eqref{Thm1} tends to the simple expression
\[\int_0^1\log\Gamma_q(x)dx=\log C_q+\frac{1}{2\log q}(2\zeta(2)+\log^2(q-1)).\]
The definition \eqref{Cqdef} of $C_q$ enables us to get a more simple identity. Since
\[\log C_q=-\frac{1}{12}\log q+\frac12\log(q-1)-\frac{\log^2(q-1)}{2\log q}+\log(q^{-1};q^{-1})_\infty,\]
the term $\frac{\log^2(q-1)}{2\log q}$ cancels and a trivial modification gives the second formula \eqref{Thm1spec} of our Theorem.

\section{The proof of the main formula}

The \eqref{gdammadef} definition of the $q$-Gamma function and some reduction gives that for any $q>1$ and $y>0$
\begin{equation}
\frac{1}{\Gamma_{q^2}\left(\frac12\log_q\frac qy\right)\Gamma_{q^2}\left(\frac12\log_qqy\right)}=\label{gprod}
\end{equation}
\[\frac{(q^{\frac12})^{1-\log_q^2y}}{(q^{-2};q^{-2})^3_\infty(q^2-1)}(q^{-2};q^{-2})_\infty(y/q;q^{-2})_\infty(1/(yq);q^{-2})_\infty\]
This product can be rewritten by Jacobi's triple product identity \cite[p. 15]{Gasper}:
\[(q^2;q^2)_\infty(qy;q^2)_\infty(q/y;q^2)_\infty=\sum_{n=-\infty}^\infty(-1)^nq^{n^2}y^n,\]
or, which is the same,
\[(q^{-2};q^{-2})_\infty(y/q;q^{-2})_\infty(1/(qy);q^{-2})_\infty=\sum_{n=-\infty}^\infty(-1)^nq^{-n^2}y^n.\]
In \eqref{gprod} we choose $y=q^{1-2x}$. Then $\frac12\log_q\frac qy=x$ and $\frac12\log_qqy=1-x$, so Jacobi's identity yields
\[\frac{1}{\Gamma_{q^2}(x)\Gamma_{q^2}(1-x)}=\frac{(q^{\frac12})^{1-(1-2x)^2}}{(q^{-2};q^{-2})^3_\infty(q^2-1)}\sum_{n=-\infty}^\infty(-1)^nq^{-n^2}(q^{1-2x})^n.\]
Consider the definition of the $\vartheta_4$ function on the first page. It is not hard to see that we arrive at the next formula:
\[\frac{1}{\Gamma_{q^2}(x)\Gamma_{q^2}(1-x)}=\frac{q^{2x(1-x)}}{(q^{-2};q^{-2})^3_\infty(q^2-1)}\vartheta_4\left(\frac{1}{2i}(1-2x)\log q,\frac1q\right).\]
In the next step we take logarithm of both sides and integrate on $[0,1]$.
\[\int_0^1\log\Gamma_{q^2}(x)\Gamma_{q^2}(1-x)dx=\]
\[\log(q^{-2};q^{-2})^3_\infty(q^2-1)-\log q\int_0^12x(1-x)dx-\int_0^1\log\vartheta_4\left(\frac{1}{2i}(1-2x)\log q,\frac1q\right)dx.\]
Using Theorem 1, the right hand side equals to
\[\frac{2\zeta(2)}{\log q^2}+\log\frac{q^2-1}{\sqrt[6]{q^2}}+\log(q^{-2};q^{-2})^2_\infty.\]
An elementary simplification implies that
\[\int_0^1\log\vartheta_4\left(\frac{1}{2i}(1-2x)\log q,\frac1q\right)dx=\log(q^{-2};q^{-2})_\infty-\frac{\zeta(2)}{\log q}.\]
We transform the integral:
\[\int_0^1\log\vartheta_4\left(\frac{1}{2i}(1-2x)\log q,\frac1q\right)dx=\frac{-i}{\log q}\int_{-\frac12i\log q}^{\frac12i\log q}\log\vartheta_4\left(x,\frac1q\right)dx.\]
Therefore
\begin{equation}
\int_{-\frac12i\log q}^{\frac12i\log q}\log\vartheta_4\left(x,\frac1q\right)dx=\frac1i\left[\zeta(2)-\log q\log(q^{-2};q^{-2})_\infty\right].\label{almost}
\end{equation}
Let us consider the endpoint of the integration:
\[\vartheta_4\left(-\frac12i\log q,\frac1q\right)=\sum_{n=-\infty}^\infty(-1)^nq^{n-n^2}.\]
It is straightforward to see that all terms cancels, so
\[\vartheta_4\left(-\frac12i\log q,\frac1q\right)=0.\]
Similarly,
\[\vartheta_4\left(\frac12i\log q,\frac1q\right)=0.\]
Let $x^*=\frac12i\log q$. To visualize the roots, we draw $\vartheta_4(ix,1/2)$:

\begin{center}
\includegraphics[scale=0.8]{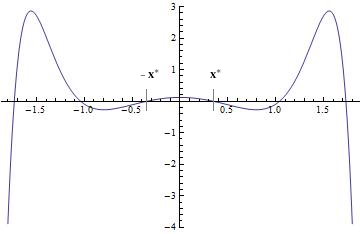}
\end{center}

Note that until this point $q>1$. Our formula will be more elegant if we substitute $q=\frac1q$. But in this case \eqref{almost} modifies:
\[\int_{x^*}^{-x^*}\log\vartheta_4\left(x,q\right)dx=\frac1i\left[\zeta(2)+\log q\log(q^2;q^2)_\infty\right].\]
Interchanging the limits of the integration and using the well known generating function \cite{Andrews}
\[\sum_{n=0}^\infty P(n)q^n=(q;q)^{-1}_\infty,\]
we are done.

\end{document}